\documentclass[11pt]{amsart}
\usepackage{enumerate}
\usepackage{amsmath, amsthm, amscd, amsfonts, amssymb, graphicx, color,float,pgf,tikz}
\usepackage{amssymb,fontenc}
\usepackage{latexsym,wasysym,mathrsfs}
\usepackage{hyperref}
\usetikzlibrary{arrows}
\usepackage[square,numbers,sort&compress]{natbib}
\makeatletter \oddsidemargin.9375in \evensidemargin \oddsidemargin
\newtheorem{theorem}{Theorem}[section]
\newtheorem{lemma}[theorem]{Lemma}
\newtheorem{proposition}[theorem]{Proposition}
\newtheorem{corollary}[theorem]{Corollary}

\theoremstyle{definition}
\newtheorem{definition}[theorem]{Definition}
\newtheorem{example}[theorem]{Example}

\newtheorem{remark}[theorem]{Remark}
\numberwithin{equation}{section}

\newcommand{\be}{\begin{equation}}
\newcommand{\ee}{\end{equation}}

\numberwithin{equation}{section}

\newcommand{\set}[1]{\{#1\}}

\newcommand{\e}{\epsilon}
\newcommand{\de}{\delta}

\newcommand{\N}{\mathbb{N}}
\newcommand{\Z}{\mathbb{Z}}
\newcommand{\R}{\mathbb{R}}
\newcommand{\m}{\mathfrak{M}}
\newcommand{\cH}{\mathcal{H}}

\newcommand{\CM}{\'Ciri\'c-Matkowski}

\usepackage{etoolbox}
\makeatletter
\patchcmd{\@settitle}{\uppercasenonmath\@title}{}{}{}
\patchcmd{\@setauthors}{\MakeUppercase}{}{}{}
\makeatother

\allowdisplaybreaks
\begin{document}


\title[Cauchy Sequences in Fuzzy Metric Spaces and \dots]%
{Cauchy Sequences in Fuzzy Metric Spaces and Fixed Point Theorems}

\author[M. Abtahi]{Mortaza Abtahi}

\address{School of Mathematics and Computer Sciences,
Damghan University, Damghan, P.O.BOX 36715-364, Iran.}

\email{abtahi\@du.ac.ir; mortaza.abtahi\@gmail.com}



\keywords{Fuzzy metric spaces, Cauchy sequences, Fixed point theorems,
Contractive mappings, Gauge functions.\\
\indent Received: dd mmmm yyyy,    Accepted: dd mmmm yyyy.
}
\maketitle
\hrule width \hsize \kern 1mm

\begin{abstract}
  In this paper, contractive mappings of \CM{} type in fuzzy metric spaces are studied.
  A class $\Psi_1$ of gauge functions $\psi:(0,1]\to(0,1]$ such that,
  for any $r\in(0,1)$, there exists $\rho\in(r,1)$ such that $1-r> \tau >1-\rho$
  implies $\psi(\tau)\geq 1-r$, is introduced, and it is shown that fuzzy $\psi$-contractive mappings
  are fuzzy contractive mappings of \CM{} type. A characterization of
  Cauchy sequences in fuzzy metric spaces is presented, and it is utilized to establish
  fixed point theorems. Examples are given to support the results. Our results
  cover those of Mihet (Fuzzy $\psi$-contractive mappings in non-Archimedean fuzzy metric spaces,
  Fuzzy Sets Syst.\ 159(2008) 739--744), Wardowski (Fuzzy contractive mappings
  and fixed points in fuzzy metric spaces, Fuzzy Sets Syst.\ 222(2013) 108--114)
  and others.
\end{abstract}
\maketitle
\vspace{0.1in}
\hrule width \hsize \kern 1mm

\section{Introduction}
\label{sec:intro}

Fuzzy metric spaces were initiated by Kramosil and Mich\'alek \cite{Kramosil-Michalek}.
Later, in order to obtain a Hausdorff topology in fuzzy metric spaces, George and Veeramani
\cite{George-Veeramani} modified the conditions formulated in \cite{Kramosil-Michalek}.
The study of fixed point theory in fuzzy metric spaces started with the work of Grabiec
\cite{Grabiec-1988}, by extending the well-known fixed point theorems of Banach \cite{Banach} and
Edelstein \cite{Edelstein} to fuzzy metric spaces. Many authors followed this concept by introducing
and investigating different types of fuzzy contractive mappings; see, e.g.,
\cite{Gregori-Sapena,Mihet-2004,Mihet-2007,Mihet-2008,Mihet-2014,Saini-2008,Sharma-2002,
Vasuki-Veeramani,Vetro-2011,Wardowski-2013,Zikic-2004}.
They reconsidered the Banach contraction principle
by initiating a new concept of fuzzy contractive mapping
in fuzzy metric spaces in the sense of George and Veeramani \cite{George-Veeramani}
and also in the sense of Kramosil and Mich\'alek \cite{Kramosil-Michalek}.
However, in their results, there were used strong conditions for completeness,
namely $G$-completeness \cite{Grabiec-1988},
of a fuzzy metric space. Being aware of this problem, they raised the question whether the fuzzy contractive
sequences are Cauchy in the usual sense, namely $M$-Cauchy \cite{George-Veeramani}.
Many papers have appeared concerning this subject; see, for example, the interesting results
of Mihet \cite{Mihet-2004,Mihet-2007,Mihet-2008}.

In \cite{Wardowski-2013} Wardowski introduced the concept of fuzzy $\cH$-contractive
mappings and formulated conditions guaranteeing the convergence of fuzzy $\cH$-contractive sequences
to a unique fixed point in a complete fuzzy metric space.
The paper includes a comprehensive set of examples showing the generality of the results and
demonstrating that the formulated conditions are significant and cannot be omitted.
However, in \cite{Gregori-Minana-2014}, it is shown that some assertions made in \cite{Wardowski-2013}
are not true.

Many works has been done on fuzzy fixed point theory and sequences of fuzzy real numbers
in the recent past; see \cite{Tripathy-1,Tripathy-2,Tripathy-3,Tripathy-4,Tripathy-5,Tripathy-6}
for example.

This paper is outlined as follows. Section \ref{sec:pre} contains definitions and basic concepts
of fuzzy metric spaces. In Section \ref{sec:psi-contractive}, we introduce and discuss contractive mappings
of \CM{} type in fuzzy metric spaces. We introduce a class $\Psi_1$
of gauge functions that contains the class $\Psi$ introduced in \cite{Mihet-2008},
and establish a relation between fuzzy \CM{} contractive mappings and fuzzy $\psi$-contractive
mappings ($\psi\in\Psi_1$). In Section \ref{sec:Cauchy}, a characterization of Cauchy sequences
in fuzzy metric spaces is presented. Utilizing this characterization, in Section \ref{sec:fp-thms},
it is proved that fuzzy \CM{} contractive mappings on complete fuzzy metric spaces
have unique fixed points. Fixed point theorems presented in Section \ref{sec:fp-thms} extend
some results on this subject (e.g., \cite{Gregori-Minana-2014,Mihet-2008,Mihet-2014,Wardowski-2013}).
At the end, to support our results, we present many examples.

\section{Preliminaries}
\label{sec:pre}

Throughout, the sets of integers, nonnegative integers and positive integers
are denoted, respectively, by $\Z$, $\Z^+$ and $\N$. The sets of real numbers
and nonnegative real numbers are denoted, respectively, by $\R$ and $\R^+$.

A binary operation $*$ on $[0,1]$ is called a \emph{triangular norm} or
a \emph{t-norm} \cite{Schweizer-Sklar}
if it is associative, commutative, and satisfies the following properties:
\begin{enumerate}
  \item $a*1=a$, for all $a\in [0,1]$,
  \item $a*b \leq c*d$ whenever $a\leq c$ and $b\leq d$.
\end{enumerate}

A t-norm $*$ is called \emph{positive} if $a*b >0$ whenever $a>0$ and $b>0$.
Some typical examples of t-norms are the following:
\begin{align*}
  a*b & = ab,   && (\text{product}), \\
  a*b & = \min\{a,b\}, &&  (\text{minimum}), \\
  a*b & = \max\{0,a+b-1\},  && (\text{Lukasiewicz}),\\
  a*b & = \frac{ab}{a+b-ab}, && (\text{Hamacher}).
\end{align*}

\begin{definition}[\cite{Kramosil-Michalek}]
   A \emph{fuzzy metric space} is a triplet $(X,M,*)$, where $X$ is a nonempty set, $*$ is
   a continuous t-norm, and $M : X^2\times [0,\infty)\to [0,1]$ is a mapping with
   the following properties:
   \begin{enumerate}[\quad(1)]
     \item \label{item:(FM1)}
     $M(x,y,0) = 0$, for all $x,y\in X$;

     \item $x=y$ if and only if $M(x,y,t) = 1$, for all $t > 0$;

     \item $M(x,y,t) = M(y,x,t)$, for all $x,y\in X$ and $t > 0$;

     \item \label{item:(FM4)}
     $M(x,y,\cdot) : [0,\infty) \to [0,1]$ is left continuous, for all $x,y \in X$;

     \item \label{item:triangle}
     $M(x,z,s+t) \geq M(x,y,s) * M(y,z,t)$, for all $x,y,z \in X$ and $s,t > 0$.
   \end{enumerate}
\end{definition}

\noindent
If, in the above definition, the triangular inequality \eqref{item:triangle} is replaced by
\begin{enumerate}[\quad($5'$)]
  \item $M(x,z,t) \geq M(x,y,t) * M(y,z,t)$, for all $x,y,z \in X$ and $t > 0$,
\end{enumerate}
then the triplet $(X,M,*)$ is called a \emph{strong fuzzy metric space}.
It is easy to check that the triangle inequality ($\ref{item:triangle}'$) implies
\eqref{item:triangle}, that is, every strong
fuzzy metric space is itself a fuzzy metric space.

In \cite{George-Veeramani}, the authors modified the above definition
in order to introduce a Hausdorff topology on the fuzzy metric space.

\begin{definition}[George and Veeramani, \cite{George-Veeramani}]
   A \emph{fuzzy metric space} is a triplet $(X,M,*)$, where $X$ is a nonempty set, $*$ is
   a continuous t-norm, and $M : X^2\times (0,\infty)\to [0,1]$ is a mapping satisfying
   the following properties, for all $x,y\in X$ and $t>0$:
   \begin{enumerate}[\quad(a)]
     \item $M(x,y,t)>0$;

     \item $x=y$ if and only if $M(x,y,t) = 1$;

     \item $M(x,y,t) = M(y,x,t)$;

     \item $M(x,y,\cdot) : (0,\infty) \to [0,1]$ is continuous;

     \item $M(x,z,s+t) \geq M(x,y,s) * M(y,z,t)$.
   \end{enumerate}
\end{definition}

For example, given a metric space $(X,d\,)$, define a t-norm by $a*b=ab$,
for all $a,b\in [0,1]$, and set
  \[
    M_d(x,y,t)=\frac t{t+d(x,y)}, \quad (x,y\in X,\, t>0).
  \]
  Then $(X,M_d,*)$ is a (strong) fuzzy metric space; $M_d$ is called the standard fuzzy
  metric induced by $d$. It is interesting to note that the topologies induced by
  the standard fuzzy metric $M_d$ and the corresponding metric $d$ coincide.
  (For more information, see \cite{George-Veeramani}.)

By definition, a sequence $(x_n)$ in a fuzzy metric space $(X,M,*)$ converges to a point $x\in X$,
if $\lim_{n\to\infty}M(x_n,x,t)=1$, for all $t>0$.
In \cite{Grabiec-1988}, notions of Cauchy sequences and complete fuzzy metric spaces
are defined. The sequence $(x_n)$ in $X$ is called \emph{$G$-Cauchy} if
\[
  \lim_{n\to\infty}M(x_n,x_{n+m},t)=1 \quad(m\in\N,t>0).
\]
A fuzzy metric space in which every $G$-Cauchy sequence is convergent is called
a \emph{$G$-complete} fuzzy metric space. With this definition of completeness, even
$(\R,M_d,*)$, where $d$ is the Euclidean metric on $\R$ and $a*b=ab$,
fails to be complete. Hence, in \cite{George-Veeramani},
the authors redefined Cauchy sequence as follows.


\begin{definition}[George and Veeramani, \cite{George-Veeramani}]\
\label{dfn:Cauchy-sequences}
  A sequence $(x_n)$ in a fuzzy metric space $(X,M,*)$ is said to be
  \emph{Cauchy} (or \emph{$M$-Cauchy}) if, for each $r\in(0,1)$ and each $t>0$,
  there exists $N\in\N$ such that
  \[
    M(x_n,x_m,t) > 1-r \quad (m,n\geq N).
  \]
  A fuzzy metric space in which every Cauchy sequence is convergent, is
  called a \emph{complete} (or \emph{$M$-complete}) fuzzy metric space.
\end{definition}

Note that a metric space $(X,d\,)$ is complete if and only if the induced
fuzzy metric space $(X,M_d,*)$ is complete.

In this paper, we work in the setting of fuzzy metric spaces in the sense of
George and Veeramani, in which by a Cauchy sequence we always mean an $M$-Cauchy sequence,
and by a complete space we always mean an $M$-complete space.

\section{Fuzzy $\psi$-contractive Mappings}
\label{sec:psi-contractive}

To begin, a result of \cite{Proinov-2006} is needed.
Let $\Phi_1$ denote the class of all functions $\phi:\R^+\to \R^+$ with the property that,
for any $\e>0$, there exists $\delta>\e$ such that $\e< s <\delta$ implies $\phi(s)\leq\e$.

\begin{lemma}[\cite{Proinov-2006}]
\label{lem:property-of-Phi1}
  Let $E$ and $F$ be two nonnegative functions on a nonempty set $X$ such that $E(x)\leq F(x)$,
  for all $x\in X$. Then the following statements are equivalent:
  \begin{enumerate}[\upshape(i)]
    \item There is a function $\phi\in\Phi_1$ such that $E(x)\leq \phi(F(x))$, for all $x\in X$.
    \item For any $\e>0$ there is $\delta>\e$ such that $\e <F(x)<\delta$ implies $E(x)\leq\e$.
  \end{enumerate}
\end{lemma}

In \cite{Mihet-2008} (see also \cite{Mihet-2004}) Mihet defined a class $\Psi$ consisting
of all continuous, nondecreasing functions $\psi: (0,1] \to (0,1]$ with the property that
$\psi(\tau) > \tau$, for all $\tau\in(0,1)$. In the following, we introduce a class $\Psi_1$
of gauge functions that includes $\Psi$.

\begin{definition}
\label{dfn:psi-contractive}
   Let $\Psi_1$ denote the class of all functions $\psi:(0,1]\to (0,1]$ with the property that,
   for any $r\in(0,1)$, there exists $\rho\in(r,1)$ such that $1-r> \tau >1-\rho$ implies $\psi(\tau)\geq 1-r$.
   Given $\psi\in\Psi_1$, a self-map $T$ of a fuzzy metric space $(X,M,*)$ is said to be
   \emph{$\psi$-contractive} if
   \begin{enumerate}
     \item $M(Tx,Ty,t)>M(x,y,t)$ for $x\neq y$ and $t>0$,
     \item $M(Tx,Ty,t)\geq \psi(M(x,y,t))$, for all $x,y$ and $t>0$.
   \end{enumerate}
\end{definition}

We prove fixed point theorems for $\psi$-contractive mappings ($\psi\in\Psi_1$)
in complete fuzzy metric spaces (Theorems \ref{thm:CM:strong} and
\ref{thm:CM:general}). Since the class $\Psi_1$ properly contains the class $\Psi$
(Proposition \ref{prop:Psi subset Psi1} and Example \ref{exa:Psi1 properly contains Psi}),
our theorem extends Mihet's result in \cite{Mihet-2008}.

\begin{proposition}\label{prop:Psi subset Psi1}
  $\Psi \subset\Psi_1$.
\end{proposition}

\begin{proof}
  Let $\psi\in\Psi$. For every $r\in(0,1)$, since $\psi(1-r)>1-r$, the continuity of $\psi$ at $1-r$
  implies the existence of some $\rho\in (r,1)$ such that $1-r> \tau >1-\rho$ implies $\psi(\tau)\geq 1-r$.
  This shows that $\psi\in\Psi_1$, and thus $\Psi\subset\Psi_1$.
\end{proof}

\begin{remark}
  In Example \ref{exa:Psi1 properly contains Psi}, a function $\psi\in\Psi_1$
  which is not continuous is presented. This confirms that the inclusion $\Psi \subset \Psi_1$ is proper.
\end{remark}

As in \cite{Wardowski-2013},
let $\cH$ be the family of all strictly decreasing bijections $\eta:(0,1]\to[0,\infty)$.
Wardowski \cite{Wardowski-2013} introduced the notion of fuzzy $\cH$-contractive mappings,
and proved a fixed point theorem for such mappings. In \cite{Gregori-Minana-2014}, it is
shown that fuzzy $\cH$-contractive mappings are included in the class of fuzzy $\psi$-contractive
mappings ($\psi\in\Psi$). One might define new fuzzy contractive mappings by composing functions
$\eta\in\cH$ with gauge functions $\phi\in\Phi_1$. Given $\eta\in\cH$ and
$\phi\in\Phi_1$, let us call a self-map $T$ of a fuzzy metric space $(X,M,*)$
a \emph{fuzzy $\phi$-$\cH$-contractive mapping}, if
  \[
    \eta(M(Tx,Ty,t)) \leq \phi\bigl( \eta(M(x,y,t)) \bigr)
    \quad (x,y\in X,t>0).
  \]

We show that these kinds of fuzzy contractive mappings are included in the class
of fuzzy $\psi$-contractive mappings ($\psi\in\Psi_1$).

\begin{proposition}
\label{prop:Psi1=Phi1}
  For every $\eta\in \cH$, the mapping $\phi\mapsto \eta^{-1}\circ \phi \circ \eta$, is
  a one-to-one correspondence from $\Phi_1$ onto $\Psi_1$.
\end{proposition}

\begin{proof}
  Let $\phi\in \Phi_1$, and $\psi=\eta^{-1}\circ \phi \circ \eta$.
  Given $r\in(0,1)$, let $\e=\eta(1-r)$. Since $\phi\in \Phi_1$, there
  is $\delta>\e$ such that $\e<s<\delta$ implies $\phi(s)\leq \e$.
  Let $\rho=1-\eta^{-1}(\delta)$. Now, if $1-r>\tau>1-\rho$ then
  $\e<\eta(\tau)<\delta$, and thus $\phi(\eta(\tau))\leq \e$. Therefore,
  $\psi(\tau)=\eta^{-1}(\phi(\eta(\tau)))\geq \eta^{-1}(\e)=1-r$.

  Similarly, if $\psi\in\Psi_1$, then $\phi=\eta\circ \psi \circ \eta^{-1}$
  belongs to $\Phi_1$.
\end{proof}

\begin{corollary}
  The class of fuzzy $\phi$-$\cH$-contractive mappings
  is included in the class of fuzzy $\psi$-contractive mappings, for $\psi\in\Psi_1$.
\end{corollary}

We conclude this section by introducing the concept of fuzzy \CM{} contractive
mappings.

\begin{definition}
\label{dfn:fuzzy-Ciric-Matkowski-contraction}
  A self-map $T$ of a fuzzy metric space $(X,M,*)$ is said to be
  a \emph{fuzzy \CM{} contractive mapping} if
  \begin{enumerate}
    \item $M(Tx,Ty,t) > M(x,y,t)$ for $x\neq y$ and $t>0$,

    \item for every $t>0$ and $r\in(0,1)$, there exists $\rho\in(r,1)$ such that
    \begin{equation}\label{eqn:fuzzy-Ciric-Matkowski-contraction}
      1-r>M(x,y,t)>1-\rho \quad \Rightarrow \quad M(Tx,Ty,t) \geq 1-r \quad(x,y\in X).
    \end{equation}
  \end{enumerate}
\end{definition}

\noindent
Note that, in the above definition, condition \eqref{eqn:fuzzy-Ciric-Matkowski-contraction}
can be replaced by the following:
\begin{equation}\label{eqn:fuzzy-Ciric-Matkowski-contraction2}
  M(x,y,t)>1-\rho \quad \Rightarrow \quad M(Tx,Ty,t) \geq 1-r \quad (x,y\in X).
\end{equation}

The following reveals the relation between fuzzy $\psi$-contractive
mappings and fuzzy \CM{} contractive mappings.

\begin{lemma}
\label{lem:property-of-Psi1}
  Let $E$ and $F$ be two nonnegative functions on a nonempty set $X$
  such that $E(x)\geq F(x)$, for all $x\in X$. Then the following statements
  are equivalent:
  \begin{enumerate}[\upshape(i)]
    \item \label{item:psi-contraction}
    There is a function $\psi\in\Psi_1$ such that $E(x)\geq \psi(F(x))$,
    for all $x\in X$.

    \item \label{item:CM-contraction}
    For any $r\in(0,1)$ there is $\rho\in(r,1)$ such that
    \begin{equation}\label{eqn:F-E}
       1-r>F(x)>1-\rho \quad \Rightarrow \quad E(x)\geq1-r \quad (x\in X).
    \end{equation}
  \end{enumerate}
\end{lemma}

\begin{proof}
  $\eqref{item:psi-contraction} \Rightarrow \eqref{item:CM-contraction}$
  Clear.

  $\eqref{item:CM-contraction} \Rightarrow \eqref{item:psi-contraction}$
  Let $\eta\in\cH$, and set $E_1(x)=\eta(E(x))$ and $F_1(x)=\eta(F(x))$.
  It is a matter of calculation to see that, for any $\e>0$, there is $\delta>\e$
  such that $\e <F_1(x)<\delta$ implies $E_1(x)\leq\e$. Hence, by Lemma \ref{lem:property-of-Phi1},
  there is $\phi\in\Phi_1$ such that $E_1(x) \leq \phi(F_1(x))$, for all $x\in X$.
  Therefore,
  \[
    E(x)=\eta^{-1}(E_1(x)) \geq \eta^{-1}(\phi(F_1(x))) = \eta^{-1}\bigl(\phi(\eta(F_1(x)))\bigr).
  \]
  If we take $\psi=\eta^{-1}\circ \phi \circ \eta$, then $\psi\in\Psi_1$
  (by Proposition \ref{prop:Psi1=Phi1}) and $E(x)\geq \psi(F(x))$,
  for all $x\in X$.
\end{proof}

Obviously, the class of fuzzy \CM{} contractive mappings contains
the class of $\psi$-contractive mappings ($\psi\in\Psi_1$). We can say more:

\begin{theorem}
\label{thm:equivalences}
  Let $T$ be a self-map of a fuzzy metric space $(X,M,*)$ such that
  \[
    M(Tx,Ty,t)\geq M(x,y,t) \quad (x,y\in X, t>0).
  \]
  Consider the following statements;
\begin{enumerate}[\upshape(i)]
  \item \label{item:T-is-psi-contractive}
  There exists $\psi\in\Psi_1$ such that
  \[
    M(Tx,Ty,t) \geq \psi(M(x,y,t)) \quad (x,y\in X,t>0).
  \]

  \item \label{item:T-is-uni-CM-contractive}
  For every $r\in(0,1)$, there exists $\rho\in (r,1)$ such that
  \[
    1-r>M(x,y,t)>1-\rho \quad \Rightarrow \quad M(Tx,Ty,t)\geq 1-r \quad (x,y\in X,t>0).
  \]

  \item \label{item:T-is-pw-CM-contractive}
  For every $t>0$ and $r\in(0,1)$, there exists $\rho\in (r,1)$ such that
  \[
    1-r>M(x,y,t)>1-\rho \quad \Rightarrow \quad M(Tx,Ty,t)\geq 1-r \quad (x,y\in X).
  \]

  \item \label{item:T-is-psi-t-contractive}
  For every $t>0$, there exists $\psi_t\in\Psi_1$ such that
  \[
    M(Tx,Ty,t) \geq \psi_t(M(x,y,t)) \quad (x,y\in X).
  \]
\end{enumerate}
Then $\eqref{item:T-is-psi-contractive} \Leftrightarrow \eqref{item:T-is-uni-CM-contractive}
  \Rightarrow \eqref{item:T-is-pw-CM-contractive} \Leftrightarrow
  \eqref{item:T-is-psi-t-contractive}$
\end{theorem}

\begin{proof}
  $\eqref{item:T-is-psi-contractive} \Leftrightarrow \eqref{item:T-is-uni-CM-contractive}$
  follows from Lemma \ref{lem:property-of-Psi1} by considering the nonnegative functions
  \[
    E(x,y,t)=M(Tx,Ty,t),\ F(x,y,t)=M(x,y,t),
  \]
  on $X^2\times(0,\infty)$. The implication
  $\eqref{item:T-is-uni-CM-contractive} \Rightarrow \eqref{item:T-is-pw-CM-contractive}$
  is obvious, and $\eqref{item:T-is-pw-CM-contractive} \Leftrightarrow \eqref{item:T-is-psi-t-contractive}$
  also follows from Lemma \ref{lem:property-of-Psi1} by considering, for each $t>0$,
  the nonnegative functions $E_t(x,y)=M(Tx,Ty,t)$ and $F_t(x,y)=M(x,y,t)$
  on $X\times X$.
\end{proof}

\section{A Characterization of Cauchy Sequences}
\label{sec:Cauchy}

In this section, a characterization of Cauchy sequences in fuzzy metric spaces
is presented. This characterization is used, in Section \ref{sec:fp-thms},
to give new fixed point theorems. First, we need a couple of 
definitions.

\begin{definition}
\label{dfn:asym-reg-seq}
  Let $(x_n)$ be a sequence in a fuzzy metric space $(X,M,*)$.
  \begin{enumerate}
    \item The sequence $(x_n)$ is called \emph{asymptotically regular} if $M(x_n,x_{n+1},t)\to1$,
    for every $t>0$.

    \item The sequence $(x_n)$ is called \emph{uniformly asymptotically regular} if,
  for any sequence $E=\set{t_i:i\in\N}$ of positive numbers with $t_i\searrow 0$,
  we have $M(x_n,x_{n+1},t) \to 1$, uniformly on $t\in E$.
  \end{enumerate}
\end{definition}

The following is the main result of the section.

\begin{lemma}
\label{lem:main-lem}
  Let $(x_n)$ be a sequence in $(X,M,*)$. Suppose,
  for every $t>0$ and $r\in(0,1)$, for any two subsequence $(x_{p_n})$ and $(x_{q_n})$,
  if\/ $\liminf M(x_{p_n},x_{q_n},t)\geq1-r$, then, for some $N$,
  \begin{equation}
    M(x_{p_n+1},x_{q_n+1},t)\geq1-r, \quad (n\geq N).
  \end{equation}
  Then $(x_n)$ is Cauchy, provided it is uniformly asymptotically regular. In case
  $(X,M,*)$ is a strong fuzzy metric space, we only need $(x_n)$ be asymptotically regular.
\end{lemma}

\begin{proof}
  To get a contradiction, assume that $(x_n)$ is not Cauchy. Then, there exist $t>0$
  and $r\in(0,1)$ such that
  \begin{equation}\label{eqn:negation-of-Cauchy}
    \forall k\in\N,\ \exists\, p,q\geq k,
    \quad M(x_p,x_q,t)<1-r.
  \end{equation}

  \noindent
  Set $a_i=1/i$, and $E=\set{a_it: i\in\N}$. In case $(x_n)$ is uniformly
  asymptotically regular, we have
  $M(x_n,x_{n+1},a_it)\to1$, uniformly on $E$. Hence,
  there exist positive integers $k_1<k_2<\dotsb$ such that
  \begin{equation}\label{eqn:uni-asy-reg}
    M(x_m,x_{m+1},a_it) > 1 - \frac rn, \quad (m\geq k_n,\, i\in \N).
  \end{equation}

  \noindent
  For each $k_n$, by \eqref{eqn:negation-of-Cauchy}, there exist
  integers $p_n$ and $q_n$ such that
  \begin{equation}\label{eqn:pn-and-qn}
    q_n>p_n\geq k_n\quad \text{and} \quad
    M(x_{p_n+1},x_{q_n+1},t) < 1-r.
  \end{equation}

  \noindent
  We choose $q_n$ be the smallest such integer so that $M(x_{p_n+1},x_{q_n},t)\geq 1-r$.
  Now, for every $n$ and $i$, we have
  \begin{align}
    M(x_{p_n},x_{q_n},t)
     & \geq M(x_{p_n},x_{p_n+1},a_it) * M(x_{p_n+1},x_{q_n},(1-a_i)t) \label{eqn:Tri-Ineq}\\
     & > \Bigl(1-\frac rn\Bigr)*M(x_{p_n+1},x_{q_n},(1-a_i)t) \nonumber.
  \end{align}
  This is true for every $i\in\N$. If $i\to\infty$ then $(1-a_i)t\to t$, and the continuity
  of $M$ gives
  \begin{align*}
    M(x_{p_n},x_{q_n},t)
      &\geq \Bigl(1-\frac rn\Bigr) * M(x_{p_n+1},x_{q_n},t) \\
      & \geq \Bigl(1-\frac rn\Bigr)*(1-r).
  \end{align*}

  \noindent
  This implies that $\liminf M(x_{p_n},x_{q_n},t) \geq 1-r$. However,
  we have
  \[
    M(x_{p_n+1},x_{q_n+1},t) < 1-r, \quad(n\in\N).
  \]
  This is a contradiction.

  In case $X$ is a strong fuzzy metric space, we choose $k_1<k_2<\dotsb$ such that
  \eqref{eqn:uni-asy-reg} holds only for $i=1$, that is,
  \begin{equation}\label{eqn:pw-asy-reg}
    M(x_m,x_{m+1},t) > 1 - \frac rn, \quad (m\geq k_n).
  \end{equation}

  \noindent
  Then, instead of \eqref{eqn:Tri-Ineq}, we have the following
  \begin{align}
    M(x_{p_n},x_{q_n},t)
     & \geq M(x_{p_n},x_{p_n+1},t) * M(x_{p_n+1},x_{q_n},t) \label{eqn:Tri-Inequ} \\
     & > \Bigl(1-\frac rn\Bigr)*M(x_{p_n+1},x_{q_n},t). \nonumber
  \end{align}

  \noindent
  The rest of the proof is similar.
\end{proof}

The following result follows directly from the above lemma.

\begin{theorem}
\label{thm:main-thm}
  Let $(x_n)$ be a sequence in $(X,M,*)$, and
  let $\m(x,y,t)$ be a nonnegative function on $X^2\times (0,\infty)$ such that,
  for any two subsequences $(x_{p_n})$ and $(x_{q_n})$,
  \begin{equation}\label{eqn:limsup m <= limsup d}
      \liminf_{n\to\infty} \m(x_{p_n},x_{q_n},t)
      \geq \liminf_{n\to\infty} M(x_{p_n},x_{q_n},t),\quad (t>0).
  \end{equation}
  Suppose, for every $t>0$ and $r\in(0,1)$, for any two subsequences
  $(x_{p_n})$ and $(x_{q_n})$, condition
  \[
    \liminf \m(x_{p_n},x_{q_n},t)\geq 1-r,
  \]
  implies that, for some $N\in\N$,
  \begin{equation}\label{eqn:M geq 1-r}
    M(x_{p_n+1},x_{q_n+1},t) \geq 1-r, \quad (n\geq N).
  \end{equation}
  Then $(x_n)$ is Cauchy, provided it is uniformly asymptotically regular.
  In case $(X,M,*)$ is a strong fuzzy metric space, we only need $(x_n)$
  be asymptotically regular.
\end{theorem}

\begin{proof}
  Using Lemma \ref{lem:main-lem}, let $t>0$ and $r\in(0,1)$, and let $(x_{p_n})$ and $(x_{q_n})$
  be two subsequences of $(x_n)$ with
  \[
     \liminf M(x_{p_n},x_{q_n},t)\geq 1-r.
  \]
  Then $\liminf \m(x_{p_n},x_{q_n},t)\geq 1-r$, and thus \eqref{eqn:M geq 1-r} holds.
  All conditions in Lemma \ref{lem:main-lem} are fulfilled and so the sequence
  is Cauchy.
\end{proof}

The following result helps us apply Lemma \ref{lem:main-lem} and Theorem \ref{thm:main-thm}
to fuzzy \CM{} contractive mappings. In fact, if $T$ is such a contraction, then
the Picard iterations $x_n=T^nx$, $n\in\N$, satisfy condition \eqref{item:it-is-m-contractive-sequence}
in the following lemma.

\begin{lemma}
\label{lem:equiv-conditions-for-m-contractive-sequences}
 Let $(x_n)$ be a sequence in $(X,M,*)$. For a nonnegative function  $\m(x,y,t)$
 on $X^2\times (0,\infty)$, the following statements are equivalent:
  \begin{enumerate}[\upshape(i)]
    \item \label{item:it-is-m-contractive-sequence}
    for every $t>0$ and $r\in(0,1)$, there exists $\rho\in(r,1)$ and $N\in\Z^+$ such that
    \begin{equation}\label{eqn:m-contractive-sequence}
      \forall p,q\geq N, \quad
      \m(x_p,x_q,t) > 1-\rho \quad \Rightarrow \quad M(x_{p+1},x_{q+1},t) \geq 1-r.
    \end{equation}

    \item \label{item:m-contractive-sequence-in-term-of-pnqn}
    for every $t>0$ and $r\in(0,1)$, for any two subsequences $(x_{p_n})$ and $(x_{q_n})$,
    if $\liminf \m(x_{p_n},x_{q_n},t)\geq 1-r$ then,
    for some $N$,
    \[
      M(x_{p_n+1},x_{q_n+1},t) \geq 1-r, \quad (n\geq N).
    \]
  \end{enumerate}
\end{lemma}

\begin{proof}
 $\eqref{item:it-is-m-contractive-sequence}
 \Rightarrow \eqref{item:m-contractive-sequence-in-term-of-pnqn}$
 Let $t>0$ and $r\in(0,1)$. Assume, for subsequences $(x_{p_n})$ and $(x_{q_n})$,
 we have $\liminf \m(x_{p_n},x_{q_n},t)\geq 1-r$.
 By \eqref{item:it-is-m-contractive-sequence}, there exists
 $\rho\in(r,1)$ and $N_1\in\Z^+$ such that \eqref{eqn:m-contractive-sequence}
 holds. Let $N_2\in\Z^+$ be such that $\m(x_{p_n},x_{q_n},t) > 1-\rho$
 for $n\geq N_2$. Then
 \[
    M(x_{p_n+1},x_{q_n+1},t) \geq 1-r, \quad (n\geq \max\set{N_1,N_2}).
 \]

 $\eqref{item:m-contractive-sequence-in-term-of-pnqn}
 \Rightarrow \eqref{item:it-is-m-contractive-sequence}$
 Assume, to get a contradiction, that \eqref{item:it-is-m-contractive-sequence}
 fails to hold. Then there exist $t>0$, $r\in(0,1)$, and subsequences $(x_{p_n})$ and $(x_{q_n})$
 such that
 \[
   \m(x_{p_n},x_{q_n},t) > (1-r)*\Bigl(1-\frac1n\Bigr) \quad \text{and} \quad
   1-r > M(x_{p_n+1},x_{q_n+1},t).
 \]
 This contradicts \eqref{item:m-contractive-sequence-in-term-of-pnqn}
 because $\liminf \m(x_{p_n},x_{q_n},t) \geq 1-r$.
\end{proof}

\section{Fixed Point Theorems}
\label{sec:fp-thms}

In this section, we use the characterization of Cauchy sequences presented in previous section
to prove new fixed point theorems for \CM{} contractive mappings on complete fuzzy metric spaces.

\begin{definition}
  A self-map $T$ of a fuzzy metric space $(X,M,*)$ is called (\emph{uniformly}) \emph{asymptotically regular}
  at $x\in X$, if the sequence $\{T^nx\}$ is (uniformly) asymptotically regular
  (in the sense of Definition \ref{dfn:asym-reg-seq}).
\end{definition}

\begin{lemma}
\label{lem:if-T-is-CM-cont-then-T-is-asym-reg}
  Every fuzzy \CM{} contractive mapping $T$ on $X$ is asymptotically regular
  at each $x\in X$.
\end{lemma}

\begin{proof}
   Let $x\in X$ and set $x_n = T^nx$, $n\in\N$. If $x_m=x_{m+1}$, for some $m$, then
   $x_m=x_n$ for all $n\geq m$, and there is nothing to prove. Suppose $x_n\neq x_{n+1}$,
   for all $n$. Then, by induction, we have $M(x_{n+1},x_{n+2},t)> M(x_n,x_{n+1},t) > 0$,
   for all $t > 0$. Therefore, for every $t > 0$, the sequence $M(x_n,x_{n+1},t)$
   converges to some number $L(t)\in(0,1]$, and $M(x_n,x_{n+1},t)<L(t)$. We show that $L(t)=1$.
   If $L(t)<1$, take $r=1-L(t)$. There is $\rho\in(r,1)$ such that
   \eqref{eqn:fuzzy-Ciric-Matkowski-contraction} holds. Therefore,
   \[
     L(t)>M(x_n,x_{n+1},t)>1-\rho \quad \Rightarrow \quad
     M(x_{n+1},x_{n+2},t) \geq L(t).
   \]
   This is a contradiction, and thus $L(t)=1$.
\end{proof}

The following two fixed point theorems follow directly from
Lemma \ref{lem:main-lem},
Lemma \ref{lem:equiv-conditions-for-m-contractive-sequences}, and
Lemma \ref{lem:if-T-is-CM-cont-then-T-is-asym-reg}.
\begin{theorem}
\label{thm:CM:strong}
  Let $(X,M,*)$ be a complete strong fuzzy metric space. Then every fuzzy \CM{} contractive
  mapping (in particular, every $\psi$-contractive mapping for $\psi\in\Psi_1$)
  on $X$ has a unique fixed point.
\end{theorem}

\begin{theorem}
\label{thm:CM:general}
  Let $(X,M,*)$ be a complete fuzzy metric space. Then every fuzzy \CM{} contractive
  mapping (in particular, every $\psi$-contractive mapping for $\psi\in\Psi_1$)
  on $X$ has a unique fixed point provided $T$ is uniformly asymptotically regular
  at some point $x_0$.
\end{theorem}

%

\noindent
  One crucial condition in the main theorem in \cite{Wardowski-2013} is the following.
  \begin{itemize}\itshape
    \item $\set{\eta(M(x,Tx,t_i)):i\in\N}$ is bounded, for all $x\in X$ and any sequence
    $(t_i)$ of positive numbers with $t_i \searrow 0$.
  \end{itemize}
  It is worth noting that this condition implies that $T$ is
  uniformly asymptotically regular at each point $x\in X$. We see that
  Theorem \ref{thm:CM:general} also generalizes \cite[Theorem 3.2]{Wardowski-2013}.

We now present our final fixed point theorem. For a self-map $T$ of $(X,M,*)$ and nonnegative
real numbers $\alpha,\beta$, define a function $\m$ on $X^2\times(0,\infty)$ by
\begin{equation}\label{eqn:m}
  \m(x,y,t)=M(x,y,t)*M(x,Tx,t)^\alpha *M(y,Ty,t)^\beta.
\end{equation}

\begin{theorem}
\label{thm:final}
  Let $(X,M,*)$ be a complete fuzzy metric space, and $T$ be a continuous self-map of
  $X$. Define $\m$ by \eqref{eqn:m}, and suppose
  \begin{enumerate}[\upshape(i)]
    \item \label{item:strict-ineq}
    $M(Tx,Ty,t)>\m(x,y,t)$, for $x\neq y$ and $t>0$;

    \item \label{item:m-contractive}
    For every $t>0$ and $r\in(0,1)$, there exist $\rho\in(r,1)$ and
    $N\in\Z^+$ such that, for every $x,y\in X$,
    \begin{equation}\label{eqn:m-contractive}
      \m(T^Nx,T^Ny,t)>1-\rho \quad \Rightarrow \quad M(T^{N+1}x,T^{N+1}y,t)\geq 1-r.
    \end{equation}
  \end{enumerate}
  Then $T$ has a unique fixed point, provided $T$ is uniformly asymptotically regular
  at some $x_0\in X$. In case $(X,M,*)$ is a strong fuzzy metric space,
  we only need $T$ be asymptotically regular at $x_0$.
\end{theorem}

\begin{proof}
  First, let us show that $T$ has at most one fixed point. Suppose $x=Tx$ and $y=Ty$.
  By \eqref{eqn:m}, we have $\m(x,y,t)=M(x,y,t)=M(Tx,Ty,t)$, for all $t>0$.
  Now, condition \eqref{item:strict-ineq} implies that $x=y$.

  Suppose $T$ is (uniformly) asymptotically regular at $x_0$, and set $x_n=T(x_{n-1})$,
  for $n\geq1$. Then $M(x_n,x_{n+1},t)\to1$, for all $t>0$. Therefore, for any two
  subsequences $(x_{p_n})$ and $(x_{q_n})$, we have
  \begin{align*}
    \liminf\m&(x_{p_n},x_{q_n},t) \\
      & = \liminf M(x_{p_n},x_{q_n},t)* M(x_{p_n},x_{p_{n+1}},t)^\alpha* M(x_{q_n},x_{q_{n+1}},t)^\beta\\
      & = \liminf M(x_{p_n},x_{q_n},t) * 1 *1 \\
      & = \liminf M(x_{p_n},x_{q_n},t),
  \end{align*}
  which means that condition
  \eqref{eqn:limsup m <= limsup d} in Theorem \ref{thm:main-thm} holds. Moreover,
  \eqref{eqn:m-contractive} implies \eqref{eqn:m-contractive-sequence} in
  Lemma \ref{lem:equiv-conditions-for-m-contractive-sequences}. All conditions in
  Theorem \ref{thm:main-thm} are fulfilled, and hence the Picard iterations
  $T^nx_0$ form a Cauchy sequence. Since $X$ is complete, there is $z\in X$ with $T^nx_0\to z$.
  Since $T$ is continuous, we get $Tz=z$.
\end{proof}

\begin{corollary}
\label{cor:final}
  Let $(X,M,*)$ be a complete fuzzy metric space, and $T$ be a continuous self-map of
  $X$. Define $\m$ by \eqref{eqn:m}, and suppose
  \begin{enumerate}[\upshape(i)]
    \item $M(Tx,Ty,t)>\m(x,y,t)$, for $x\neq y$ and $t>0$;

    \item for every $t>0$, there exists $\psi_t\in\Psi_1$ such that, for every $x,y\in X$,
    \begin{equation*}
      M(Tx,Ty,t)\geq \psi_t(\m(x,y,t)).
    \end{equation*}
  \end{enumerate}
  Then $T$ has a unique fixed point, provided $T$ is uniformly asymptotically regular
  at some $x_0\in X$. In case $(X,M,*)$ is a strong fuzzy metric space,
  we only need $T$ be asymptotically regular at $x_0$.
\end{corollary}

\begin{corollary}
\label{cor:finall}
  Let $(X,M,*)$ be a complete fuzzy metric space, and $T$ be a continuous self-map of
  $X$. Define $\m$ by \eqref{eqn:m}, and suppose there exists $\psi\in\Psi$ such that
  \begin{equation*}
      M(Tx,Ty,t)\geq \psi(\m(x,y,t)) \quad (x,y\in X,t>0).
  \end{equation*}
  Then $T$ has a unique fixed point, provided $T$ is uniformly asymptotically regular
  at some $x_0\in X$. In case $(X,M,*)$ is a strong fuzzy metric space,
  we only need $T$ be asymptotically regular at $x_0$.
\end{corollary}

\bigskip

To justify our results, we now present many examples.

\begin{example}
\label{exa:Psi1 properly contains Psi}
  This example stems from \cite[Example 1]{Jachymski-1995}. We present a gauge
  function $\psi\in\Psi_1$ that is not continuous. Beside Proposition \ref{prop:Psi subset Psi1},
  this shows that $\Psi_1$ properly contains $\Psi$.
  Define $\psi:(0,1]\to(0,1]$ by
  \begin{equation}\label{psi(tau)}
    \psi(\tau)=
    \left\{\!\!
      \begin{array}{ll}
        \frac12, & \hbox{$\tau<\frac12$};\\[1ex]
        \frac{n+1}{n+2}, & \hbox{$\frac{n}{n+1}\leq\tau<\frac{n+1}{n+2}$;} \\[1ex]
        1, & \hbox{$\tau=1$.}
      \end{array}
    \right.
  \end{equation}

  Obviously, $\psi$ is not continuous and thus $\psi\notin\Psi$. To show that $\psi\in\Psi_1$,
  first we see that $\psi$ is nondecreasing and satisfies the following conditions:
  \begin{equation}
    \psi(\tau)>\tau\quad \text{and} \quad \lim_{n\to\infty}\psi^n(\tau)=1, \ (0<\tau<1).
  \end{equation}

  Towards a contradiction, suppose $\psi\notin\Psi_1$. Then, there exist $r\in(0,1)$ and
  a sequence $\{\tau_n\}$ in $(0,1-r)$ such that $\tau_n\to 1-r$ and $\psi(\tau_n)<1-r$.
  Let $\tau\in(0,1-r)$. Then $\tau<\tau_n<1-r$, for some $\tau_n$.
  Since $\psi$ is nondecreasing, we get $\psi(\tau)\leq \psi(\tau_n)<1-r$.
  Replacing $\tau$ by $\psi(\tau)$, we get $\psi^2(\tau)<1-r$. By induction, we have
  $\psi^n(\tau)<1-r$, for all $\tau\in(0,1-r)$ and $n\in\N$. This contradicts the fact that
  $\psi^n(\tau)\to1$, for all $\tau\in(0,1)$, and so $\psi\in\Psi_1$.

   In fact, $\psi=\eta^{-1}\circ \phi\circ\eta$,
   where $\eta(\tau)=1/\tau-1$, and $\phi\in\Phi_1$ is given by
  \begin{equation}\label{eqn:phi(s)}
    \phi(s)=
    \left\{\!\!
      \begin{array}{ll}
        0, & \hbox{$s=0$;} \\[1ex]
        \frac1{n+1}, & \hbox{$\frac1{n+1}<s\leq\frac1n$;} \\[1ex]
        1, & \hbox{$s>1$.}
      \end{array}
    \right.
  \end{equation}
\end{example}

The following example shows that Theorems \ref{thm:CM:strong} and \ref{thm:CM:general}
are genuine extensions of \cite[Theorem 3.1]{Mihet-2008}.

\begin{example}
\label{exa:our-thm-extends-Mihet}
  Let $X=[0,\infty)$ and define $d$ as follows \cite{Jachymski-1995}
  \begin{equation*}
    d(x,y)=
    \left\{\!\!
      \begin{array}{ll}
        \max\{x,y\}, & \hbox{$x\neq y$;} \\
        0, & \hbox{$x=y$.}
      \end{array}
    \right.
  \end{equation*}

  Then $(X,d\,)$ is a complete metric space. Let $*$ be the product
  t-norm, and $M_d$ be the fuzzy metric on $X$ induced by $d$; that is
  \[
    M_d(x,y,t)=\frac t{t+d(x,y)}, \quad (x,y\in X,\, t>0).
  \]
  Then $(X,M_d,*)$ is a complete strong fuzzy metric space.
  Define $\phi:X\to X$ as in \eqref{eqn:phi(s)}.
  Then $\phi$ is a continuous self-map of the fuzzy metric space $(X,M_d,*)$. In fact, if
  $M_d(x_n,x,t)\to1$, for $t>0$, then $d(x_n,x)\to0$ and thus $x=0$ and $x_n\to0$. It is then
  obvious that $\phi(x_n)\to0$ in $(X,M_d,*)$ which shows that $\phi$ is continuous.

  We show that $\phi$ is a fuzzy \CM{} contractive mapping. Using Theorem \ref{thm:equivalences},
  it suffices to show that there exists, for every $t>0$, a function $\psi_t\in\Psi_1$ such that
  \[
    M_d(\phi(x),\phi(y),t) \geq \psi_t(M_d(x,y,t)) \quad (x,y\in X).
  \]
  For $t>0$, define $\eta_t(\tau)=t/\tau-t$ and $\psi_t=\eta_t^{-1}\circ \phi \circ \eta_t$.
  Since $\phi\in\Phi_1$ and $\eta_t\in\cH$, Proposition \ref{prop:Psi1=Phi1} shows that
  $\psi_t\in\Psi_1$. Now, a simple calculation shows that the following statements are
  equivalent:
  \begin{enumerate}
    \item $M_d(\phi(x),\phi(y),t)\geq \psi_t(M_d(x,y,t))$, for all $x,y\in X$ and $t>0$,
    \item $d(\phi(x),\phi(y))\leq \phi(d(x,y))$, for all $x,y\in X$,
  \end{enumerate}
  and the latter is obvious. Hence, $\phi$ is a fuzzy \CM{} contractive mapping.

  Now, we show that $\phi$ fails to be a $\psi$-contractive mapping for any $\psi\in\Psi$.
  To get a contradiction, assume $\phi$ is a $\psi$-contractive mapping for some $\psi\in\Psi$.
  Since $\psi$ is continuous and $\psi(1/2)>1/2$, there exists $\rho\in(0,1/2)$ such that
  \[
    \rho < \tau < 1/2 \quad \Rightarrow \quad \psi(\tau)>1/2.
  \]
  Let $\de>0$ be such that $\rho=1/(2+\de)$. Then take $x=1$, $y=1+\de/2$ and $t=1$.
  We have
  \[
   \begin{split}
     M_d(x,y,1) & =\frac1{1+d(1,1+\de/2)} \\
                & =\frac1{2+\de/2}, \\
     M_d(\phi(x),\phi(y),1) & =\frac1{1+d(1/2,1)} \\
                & =\frac12.
   \end{split}
  \]
  We see that $\rho < M_d(x,y,1) < 1/2$ and thus we must have
  \begin{align*}
    \frac12 & = M_d(\phi(x),\phi(y),1) \\
            & \geq \psi(M_d(x,y,1)) \\
            & > \frac12,
  \end{align*}
  which is absurd.
\end{example}

\begin{example}
\label{exa:final}
  This example shows that Theorem \ref{thm:final} is a real extension of
  Theorems \ref{thm:CM:strong} and \ref{thm:CM:general}.
  Let $X=\set{0,1,2,5}$, and define $T:X\to X$ as follows:
  \[
    T0=0,\qquad T1=5,\qquad T2=0,\qquad T5=2.
  \]
  Let $*$ be the product t-norm, and set
  \[
    M(x,y,t)=\exp\Bigl[-\frac{|x-y|}t\Bigr].
  \]

  Then $(X,M,*)$ is a complete strong fuzzy metric space. It is easy to see that
  $T$ is continuous and asymptotically regular at each point of $X$. Define
  \[
    \m(x,y,t) = M(x,y,t)*M(x,Tx,t)^2*M(y,Ty,t)^2.
  \]
  Then $M(Tx,Ty,t) \geq \m(x,y,t)^{5/7}$; that is $M(Tx,Ty,t) \geq \psi(\m(x,y,t))$,
  where $\psi(\tau)=\tau^{5/7}$. Hence, all conditions in Theorem \ref{thm:final}
  (Corollary \ref{cor:finall})
  are fulfilled. However, for $x=0$ and $y=1$, we see that $M(Tx,Ty,t) < M(x,y,t)$, for all $t>0$,
  and thus Theorems \ref{thm:CM:strong} and \ref{thm:CM:general} and, in particular,
  Theorem 3.1 in \cite{Mihet-2008} do not apply to $T$.
\end{example}

\section{Conclusion}

Fuzzy metric spaces were initiated by Kramosil and Mich\'alek in \cite{Kramosil-Michalek},
and by George and Veeramani in \cite{George-Veeramani}, and some fixed point results
were obtained in \cite{Grabiec-1988}. After that, several researchers proved various
fixed point results in these spaces by investigating different types of fuzzy contractive
mappings. In the present paper, we have extended some of these results,
particularly those appearing in \cite{Wardowski-2013} and \cite{Gregori-Minana-2014}.
It has been shown, by examples, that our results are more powerful than some of the results
from the papers  \cite{Wardowski-2013} and \cite{Gregori-Minana-2014}.

\bigskip
\noindent
{\bf Acknowledgment.}
The author would like to express his sincere thanks to the referees for their valuable comments
and suggestions.

\bibliographystyle{amsplain}

\vspace{0.1in}
\hrule width \hsize \kern 1mm
\end{document}